\documentclass[11pt]{article}

\usepackage{amssymb,amsmath}
\usepackage[all,2cell]{xy}
\UseAllTwocells

\newtheorem{definition}{Definition}
\newtheorem{theorem}{Theorem}
\newtheorem{proposition}{Proposition}
\newtheorem{lemma}{Lemma}

\newenvironment{proof}{{\bf Proof. }}{\hfill $\triangleleft$}

\newcommand \id[1] {\mathrm{id}_{#1}} 
\newcommand \ot \leftarrow 
\newcommand \OT \Leftarrow 
\newcommand \TO[1] {\to_{#1}} 
\newcommand \para \parallel 
\newcommand \seq[1] {\mathrel \triangleright_{#1}} 
\newcommand \sce[1] {\sigma_{#1}} 
\newcommand \tge[1] {\tau_{#1}} 
\newcommand \SCE[2] {\sce{#1,#2}} 
\newcommand \TGE[2] {\tge{#1,#2}} 
\newcommand \unit[1]  {1_{#1}} 
\newcommand \UNIT[2]  {1_{#1,#2}} 
\newcommand \comp[1] {\ast_{#1}} 
\newcommand \free[1] {#1^{\ast}} 
\newcommand \restr[2] {#2|_{#1}}
\newcommand \extend[2] {#2^{(#1)}}
\newcommand \ZZ {\mathbb Z} 
\newcommand \NN {\mathbb N} 
\newcommand \cobnd[1] {\mathrm{i}_{#1}} 
\newcommand \HH[2] {\mathrm H_{#1}(#2)} 
\newcommand \set[1] {\{#1\}} 
\newcommand \setbis[2] {\{#1\; |\; #2\}} 
\newcommand \ctg[1] {\hbox{\bf #1}} 
\newcommand \compl {\ctg{Compl}} 
\newcommand \polyg {\ctg{Pol}}   
\newcommand \globset  {\ctg{Glob}}  
\newcommand \sets {\ctg{Sets}} 
\newcommand \glob {\ctg{O}}
\newcommand \fcompl {\ctg{Fcompl}} 
\newcommand \adjd[2] {\ar@<1ex>[d]^{#1}_{\dashv}\ar@<-1ex>@{<-}[d]_{#2}} 
\newcommand \adjr[2] {\ar@<1ex>[r]^{#1}_{\top}\ar@<-1ex>@{<-}[r]_{#2}}  
\newcommand \Hom[3] {#1(#2,#3)} 
\newcommand \opp[1] {{#1}^{\mathrm{op}}}
\newcommand \cosce[1] {{\mathrm s}_{#1}} 
\newcommand \cotge[1] {{\mathrm t}_{#1}} 
\newcommand \COSCE[2] {\cosce{#1,#2}} 
\newcommand \COTGE[2] {\cotge{#1,#2}} 
\newcommand \pair[2] {\left(#1,#2\right)} 
\newcommand \globst[1] {O[#1]} 
\newcommand \dglobst[1] {\partial{\globst{#1}}} 
\newcommand \cofib {{\mathcal I}} 
\newcommand \initial {0} 
\newcommand \abs[1]{\left|#1\right|} 
\newcommand \mclass[1] {{\mathcal {#1}}} 
\newcommand \rortho[1] {{#1}^{\perp}} 
\newcommand \lortho[1] {{}^{\perp}#1} 
\newcommand \ctxt[1] {[#1]} 
\newcommand \tctxt[2] {\ctxt{#1}_{#2}} 
\newcommand \vcx[1] {\mathbf #1 } 
\newcommand \subst[1] {{\mathrm {sub}}_#1} 
\newcommand \paracell[2] {{\mathcal{P}}_{#1}(#2)}
\newcommand \wht[2]{{\mathrm w}_{#1}(#2)} 
\newcommand \whtt[1]{\wht{}{#1}} 
\newcommand \thk[1]{\mathrm{th}(#1)} 
\newcommand \size[1]{\mathrm{size}(#1)} 
\newcommand \case[1] {\noindent $\diamond$ {\em Case~{#1}.}}
\newcommand \step[1] {\noindent $\triangleright$ {\em Step~{#1}.}}

\title{Cofibrant complexes are free}
\author{Fran\c cois M\'etayer
        \thanks{\'Equipe PPS Universit\'e Paris 7 - CNRS}}

\begin{document}

\maketitle

\abstract{We define a notion of cofibration
among $\infty$-categories and show that the cofibrant objects
are exactly the free ones, that is those generated by polygraphs.}

\section{Introduction}\label{sec:introd}

Polygraphs~\cite{burroni:higdwp,burroni:highdw}, also known as
computads~\cite{street:limicf,power:ncatpt} are structured 
systems of generators for $\infty$-categories, extending the familiar
notion of presentation by generators and relations beyond
monoids or groups, and have recently proved extremely well-adapted 
to higher-di\-men\-sional
rewriting~\cite{guiraud:trdimp,guiraud:twoppp}.

They also lead to
a simple definition of a homology
for $\infty$-categories~\cite{metayer:respol,lafontmetayer:polrhm}, based
on the following construction: 
a {\em polygraphic resolution} of an $\infty$-category
$C$ is a pair $\pair{S}{p}$ where 
\begin{itemize}
  \item $S$ is  a polygraph, generating a free $\infty$-category $\free S$;
  \item the morphism $p:\free S\to C$ is a trivial fibration
        (see~\ref{subsec:fibrat}).
\end{itemize}
$S$ gives rise to an abelian complex $\ZZ S$, whose homology only
depends on $C$, so that we may define a polygraphic homology by
\begin{displaymath}
  {\mathrm{H}}_{*}^{\mathrm{pol}}(C)=_{\mathrm{def}} \HH{*}{\ZZ S}.
\end{displaymath}
Here the main property of free $\infty$-categories is that they are 
{\em cofibrant}. In other words, given a polygraph $S$ and a
trivial fibration $p:D\to C$,
any morphism $f:\free S\to C$ lifts to a morphism
$g:\free S\to C$~(figure~\ref{fig:cofib}).
\begin{figure}[ht]
\centering
\begin{displaymath}
  \begin{xy}
    \xymatrix{ & D\ar[d]^p \\
             \free S\ar[r]_f\ar[ur]^g &C}
  \end{xy}
\end{displaymath}
\caption{lifting}
\label{fig:cofib}
\end{figure}

The purpose of the present work is to prove the converse, namely that
all cofibrant $\infty$-categories are freely generated by polygraphs,
thus establishing a simple, abstract characterization of the free objects,
otherwise defined by a rather complex inductive construction.

We first give a brief review of the basic
categories in play~(section~\ref{sec:bascat}):
$\globset$, $\compl$ and $\polyg$ stand respectively for the
category of globular sets, $\infty$-categories (or ``complexes'')
and polygraphs.  Then we investigate 
trivial fibrations and cofibrations~(section~\ref{sec:twoclasses}).
In section~\ref{sec:main}, we reduce our theorem to the fact
that the full subcategory of $\compl$ whose objects are free
is Cauchy-complete, in other words that all its idempotents split.
This is proved in appendix~\ref{annex:cauchy}. 

Let us sketch the Cauchy-completeness argument
in the simpler case of monoids: thus, let $\ctg{Mon}$ 
denote the category of monoids, and $\ctg{Fmon}$ the full subcategory
of $\ctg{Mon}$ whose objects are the free monoids. It is well-known that
a submonoid of a free monoid is not necessarily free itself. However,
if $M=\free S$ is the free monoid on the alphabet $S$ and $h:M\to M$ is an
{\em idempotent} endomorphism of $M$, then the submonoid 
${\mathrm{Fix}}(h) =\setbis{m\in M}{h(m)=m}$ of fixpoints of $h$ {\em is} free,
which easily leads to a splitting of $h$ in $\ctg{Fmon}$, hence
to the fact that $\ctg{Fmon}$ is Cauchy-complete. Here the keypoint
is to find a set of generators of ${\mathrm{Fix}}(h)$ without non-trivial
relations in $M$. A simple way to build such a set is by considering
the subset $S_1\subset S$ of those $s\in S$ such that
$h(s)=usv$ where $h(u)=h(v)=1$. Then we define $T=\setbis{h(s)}{s\in S_1}$.
It turns out that the obvious inclusion $\free T\to M$ sends $\free T$
isomorphically to ${\mathrm{Fix}}(h)$, as shown by the existence of a 
a retraction $M\to \free T$. 

Now the same ideas carry into higher
dimensions, with $\infty$-categories
instead of monoids and polygraphs instead of generating sets. 
The general case involves additional technicalities,
due to the presence of higher dimensional compositions. Of particular
importance is the notion of {\em context}, defined and explored in
appendix~\ref{annex:contexts}.   

\medskip

Many thanks to  Albert Burroni, Yves Lafont and Krzysztof Worytkiewicz,
who have been a great help in the preparation of the present work.

\section{Basic categories}\label{sec:bascat}

\subsection{Globular sets}\label{subsec:globset}

Let $\glob$ be the small category defined as follows:
\begin{itemize}
  \item the objects of $\glob$ are integers $0,1,\ldots$;
  \item the arrows are generated by composition of
        $\cosce{n},\cotge{n}:n\to n+1$, $n\in\NN$ subject to the following
        equations
        \begin{eqnarray*}
          \cosce{n+1}\circ\cosce{n} & = & \cotge{n+1}\circ\cosce{n}, \\
          \cosce{n+1}\circ\cotge{n} & = & \cotge{n+1}\circ\cotge{n}.
        \end{eqnarray*}
\end{itemize}
As a consequence, $\Hom{\glob}{m}{n}$ has exactly two elements if $m<n$, namely
$\COSCE{m}{n}=\cosce{n-1}\circ\cdots\circ\cosce{m}$ and 
$\COTGE{m}{n}=\cotge{n-1}\circ\cdots\circ\cotge{m}$.
$\Hom{\glob}{m}{n}=\emptyset$ if $m>n$, and contains the unique element
$\id{m}$ if $m=n$.
\begin{definition}\label{defin:globset}
  A {\em globular set} is a presheaf on $\glob$.
\end{definition}
In other words, a globular set is a functor from
$\opp{\glob}$ to $\sets$. Globular sets and natural transformations
form a category $\globset$. The Yoneda embedding
\begin{displaymath}
  \glob\to\globset
\end{displaymath}
takes each integer $n$ to the {\em standard globe} $\globst{n}$. We still
denote by $\cosce{n},\cotge{n}:\globst{n}\to\globst{n+1}$ 
the morphisms of globular sets representing the corresponding arrows
from $n$ to $n+1$.

Let $X$ be a globular set and $p$ an integer, the set $X(p)$ will be 
denoted by $X_p$, and its elements called
{\em cells of dimension $p$} or {\em $p$-cells}.
Hence $\globst{n}$ has exactly two $p$-cells for $p<n$, exactly one $n$-cell,
and no $p$-cells for $p>n$. Let $\dglobst{n}$ be the globular set with the
same cells as $\globst{n}$  except for $(\dglobst{n})_n=\emptyset$, and
\begin{displaymath}
  \cobnd{n}:\dglobst{n}\to\globst{n}
\end{displaymath}
the canonical injection: $\dglobst{n}$ has two $p$-cells
for $p<n$ and no other cells. 

Let us point out a few facts about $\cobnd{n}$:
\begin{itemize}
  \item $\cosce{n}\circ\cobnd{n}=\cotge{n}\circ\cobnd{n}$;
  \item there are unique maps
        $\widehat{\cosce{n}}$ and $\widehat{\cotge{n}}$ such that
        $\cosce{n}=\cobnd{n+1}\circ\widehat{\cosce{n}}$
        and $\cotge{n}=\cobnd{n+1}\circ\widehat{\cotge{n}}$;
  \item the following diagram is a pushout:
       \begin{displaymath}
      \begin{xy}
        \xymatrix{
         \dglobst{n}\ar[r]^{\cobnd{n}}\ar[d]_{\cobnd{n}} & 
          \globst{n}\ar[d]^{\widehat{\cosce{n}}} \\
         \globst{n}\ar[r]_{\widehat{\cotge{n}}}  & \dglobst{n+1} 
         }
      \end{xy}
    \end{displaymath}
\end{itemize}

Now let $X$ be a globular set, Yoneda's lemma yields a natural equivalence
\begin{equation}
  X_n \cong \Hom{\globset}{\globst{n}}{X}
 \label{eq:yoneda}
\end{equation}
hence in particular $\cosce{n}$, $\cotge{n}$ give rise to a double sequence of maps
\begin{displaymath}
  \sce{n},\tge{n}:X_{n}\OT X_{n+1}
\end{displaymath}
satisfying the {\em boundary conditions}:
\begin{eqnarray*}
  \sce{n}\circ\sce{n+1} & = & \sce{n}\circ\tge{n+1},\\
  \tge{n}\circ\sce{n+1} & = & \tge{n}\circ\tge{n+1}.
\end{eqnarray*}
Whenever $m<n$, we set $\SCE{m}{n}=\sce{m}\circ\cdots\circ\sce{n-1}$
and $\TGE{m}{n}=\tge{m}\circ\cdots\circ\tge{n-1}$. 
Let $0\leq i<n$, we say that the $n$-cells $x,y\in X_n$ are {\em $i$-composable}
if $\TGE{i}{n} x = \SCE{i}{n} y$, a relation we denote by $x\seq{i} y$.

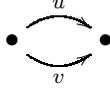
\begin{figure}[ht]
        \centering
        \begin{displaymath}
\begin{xy}
\xymatrix{
\bullet\ar @/^2ex/[r]^{u}\ar @/_2ex/[r]_{v}
&
\bullet }
\end{xy}        
\end{displaymath}
        \caption{parallel cells}
        \label{fig:parallel}
\end{figure}

If $u\in X_n$ and $\sce{n-1}(u)=x$, $\tge{n-1}(u)=y$, $x$ and $y$ are
respectively the {\em source} and the {\em target} of $u$, which we simply denote
by $u:x\to y$. Likewise, if $\SCE{i}{n}u=x$ and $\TGE{i}{n}u=y$, we shall
write $u:x\TO i y$. 
In case $u:x\to y$ and $v:x\to y$, we say that $u$, $v$ are {\em parallel},
which we denote by $u\para v$ (see figure~\ref{fig:parallel}).
 Any two $0$-cells are also considered 
to be parallel. Let $\paracell{n}{X}$ denote the set of ordered pairs
of parallel $n$-cells in $X$. We get a natural equivalence
\begin{equation}
\paracell{n}{X}\cong \Hom{\globset}{\dglobst{n+1}}{X}
 \label{eq:pairs}
\end{equation}
similar to~(\ref{eq:yoneda}). The equivalences~(\ref{eq:yoneda})
and~(\ref{eq:pairs}) assert that, for each $n$, the functors
$X\mapsto X_n$ and $X\mapsto\paracell{n}{X}$ 
from $\globset$ to $\sets$ are representable,
the representing objects being respectively $\globst{n}$ and
$\dglobst{n+1}$.

\subsection{Complexes}\label{subsec:complex}

Recall that an $\infty$-category is a globular set $C$ endowed with
\begin{itemize}
\item a {\em product} $u \comp{n-1} v : x \to z$ defined for all
$u : x \to y$ and $v : y \to z$ in $C_{n}$;
\item a {\em product} $u \comp i v : x \comp i y \to z \comp i t$
defined for all $u : x \to z$ and $v : y \to t$ in $C_{n}$ with $i < n-1$
and $u \seq i v$;
\item a {\em unit} $\unit{n{+}1} (x) : x \to x$ defined for all $x \in C_n$.
\end{itemize}
These operations satisfy the conditions of {\em associativity},
{\em left and right unit}, and {\em exchange}:
\begin{itemize}
\item$(x \comp i y) \comp i z = x \comp i (y \comp i z)$
for all $x \seq i y \seq i z$ in $C_n$ with $i < n$;
\item $\UNIT n i (x) \comp i u = u = u \comp i \UNIT n i (y)$ for all
$u : x \TO i y$ in $C_n$ with $i < n$, where
$\UNIT n i = \unit{n} \circ \unit{n{-}2} \circ \cdots \circ \unit{i+1}$;
\item $(x \comp i  y) \comp j (z \comp i t) = (x \comp j z) \comp i (y \comp j t)$
for all $x, y, z, t \in C_n$ with $i < j < n$ and $x \seq i y$, $x \seq j z$,
$y \seq j t$.
\end{itemize}
Throughout this work, {\em complex} means $\infty$-category.
Let $C$, $D$ be complexes. A {\em morphism} $f:C\to D$ is a morphism
of the underlying globular sets preserving units and products.
Complexes and morphisms build a category $\compl$, and there is
an obvious forgetful functor $\compl\to\globset$.
Its left adjoint $\globset\to\compl$ associates
to each globular set the {\em free complex} generated by it.
Note that $\globset$ is a topos of presheaves
and that the forgetful functor $\compl\to\globset$ is finitary
monadic over $\globset$. Hence $\compl$ is complete and cocomplete,
and we shall take limits and colimits in $\compl$ without
further explanations (see also~\cite{batanin:comfmg,street:pettgs}).

By restricting a complex $C$ to its cells of dimension $\leq n$, we get
an {\em $n$-category} 
\begin{displaymath}
  \restr{n}{C}:C_0\OT C_1\OT\cdots\OT C_n.
\end{displaymath}
This $n$-category can be extended to a complex $\extend{n}{C}$
by adjoining units to $\restr{n}{C}$ in all dimensions $>n$:
\begin{displaymath}
\extend{n}{C}:C_0\OT\cdots\OT C_n\OT C_n\OT\cdots.
\end{displaymath}
Let us call $\extend{n}{C}$ the {\em $n$-skeleton} of $C$.
It will be convenient to define $\extend{-1}{C}$ as the initial
complex $\initial$ with no cells.
There is a canonical inclusion
\begin{displaymath}
j_n:\extend{n}{C}\to\extend{n+1}{C}.
\end{displaymath}
Here again $j_{-1}$ denotes the unique morphism $0\to\extend{0}{C}$.

The following result is an easy consequence of the definitions:
\begin{lemma}\label{lemma:colimit}
Any complex $C$ is the colimit of its $n$-skeleta:
\begin{displaymath}
\begin{xy}
\xymatrix{\extend{-1}{C}\ar[r]^{j_{-1}} & \extend{0}{C}
\ar[r]^{j_0} & \cdots\ar[r]^{j_{n-1}} & \extend{n}{C}\ar[r]^{j_n} & \cdots}.
\end{xy}
\end{displaymath}
\end{lemma}

\subsection{Polygraphs}\label{subsec:polyg}
Let us describe a process of attaching $n{+}1$-cells
to an $n$-category $C_0\OT C_1\OT\cdots\OT C_n$.
Let $S_{n+1}$ be a set, and $\sce{n},\tge{n}:C_n\OT S_{n+1}$ a graph where
$\sce{n}$, $\tge{n}$ satisfy
the boundary conditions $\sce{n-1}\circ\sce{n}=\sce{n-1}\circ\tge{n}$
and $\tge{n-1}\circ\sce{n}=\tge{n-1}\circ\tge{n}$. We build the 
free $n{+}1$-category $C_0\OT C_1\OT\cdots\OT C_n\OT\free S_{n+1}$, where
$\free{S}_{n+1}$ consists of formal compositions of elements of
$S_{n+1}$, including identities on cells of $C_n$, and subject to the
equations of units, associativity and exchange. We refer to~\cite{burroni:highdw}
or~\cite{metayer:respol} for formal definitions.

Now $n$-polygraphs and free generated $n$-categories
are defined by simultaneous induction on $n$:
\begin{itemize}
  \item a $0$-polygraph is a set $S_0$,
        generating the  $0$-category (i.e.\ set) $\free S_0=S_0$;
  \item given an $n$-polygraph $\free S_0\OT S_1,\ldots,\free S_{n-1}\OT S_n$
        with the free $n$-category
        $\free S_0\OT\ldots\OT \free S_n$ it generates, an $n{+}1$-polygraph
        is determined by a graph $\sce{n},\tge{n}:\free S_n\OT S_{n+1}$ satisfying
        the boundary conditions, and the free $n{+}1$-category generated by it
        is $\free S_0\OT\free S_1\OT\cdots\free S_n\OT\free S_{n+1}$.  
\end{itemize}

In particular, a $1$-polygraph is simply a graph, and the notion
of free $1$-category generated by it coincides with the usual
notion of a free category generated by a graph.

\begin{definition}\label{def:polygraph}
  A {\em polygraph} is an infinite sequence
\begin{displaymath}
S:\free S_0 \OT S_1, \free S_1 \OT S_2, \ldots, \free S_n \OT S_{n{+}1}, 
\ldots
\end{displaymath}
whose first $n$ items define an $n$-polygraph for each $n$.
\end{definition}
A {\em free complex} is a complex generated by a polygraph, that is of the form
\begin{displaymath}
\free S : \free S_0  \OT  \free S_1\OT\cdots \free S_n 
\OT \free S_{n{+}1} \cdots.
\end{displaymath}

Let $S$, $T$ be polygraphs. A morphism $f:S\to T$ amounts to a sequence of maps
$f_n:S_n\to T_n$ such that, for all $\xi:x\to y$ in $S_n$,
$f_n(\xi):\free f_{n{-}1}(x)\to\free f_{n{-}1}(y)$, where $\free f_n$ is the unique extension
of $f_n$ which is compatible with products and units.

We denote by $\polyg$ the category of polygraphs and morphisms. The functor
\begin{eqnarray*}
  \polyg & \to & \compl, \\
  S & \mapsto & \free S,
\end{eqnarray*}
is left-adjoint to a forgetful functor
\begin{eqnarray*}
\compl & \to & \polyg, \\
 C & \mapsto & \abs{C}.
\end{eqnarray*}
A detailed description of $C\mapsto\abs{C}$
is given in~\cite{metayer:respol},
where this functor is called $P$.

Remark that  any globular set $X$ can be viewed as a particular polygraph
and that this identification makes $\globset$
a full subcategory of $\polyg$. Moreover the free complex generated by 
a globular set is the same as the free complex generated by the corresponding
polygraph. However most free complexes
generated by polygraphs cannot be generated by globular sets alone.

For instance the globular sets $\globst{n}$ and $\dglobst{n}$ can be 
viewed as polygraphs, and generate complexes $\free{\globst{n}}$ and
$\free{\dglobst{n}}$. Remark that in this case, the free construction
does not create new non-trivial cells. Therefore, from now on, we drop the
``$\free{}$'' in the notation of these complexes. Likewise, $\cobnd{n}$
will denote a morphism of globular sets, polygraphs, or complexes
according to the context. Note also that the natural 
equivalences~(\ref{eq:yoneda}) and~(\ref{eq:pairs}) extend
to $\compl$:
\begin{eqnarray}
C_n & \cong &\Hom{\compl}{\globst{n}}{C}
 \label{eq:cyoneda}\\
\paracell{n}{C} & \cong & \Hom{\compl}{\dglobst{n+1}}{C}
 \label{eq:cpairs} 
\end{eqnarray}
Let $S$ be a polygraph, $\free S$ the free complex it generates,
and $n$ an integer. By $\sum_{S_n} \dglobst{n}$
(resp.\ $\sum_{S_n} \globst{n}$), we mean
the direct sum of copies of $\dglobst{n}$ (resp.\ $\globst{n}$)
indexed by the elements of $S_n$. As a consequence of~(\ref{eq:cpairs}),
the source and target maps $\free S_{n-1}\OT S_n$ determine
a morphism 
\begin{displaymath}
 \rho:\sum_{S_n} \dglobst{n}\to \extend{n-1}{(\free S)}. 
\end{displaymath}
Then the following result is merely a reformulation of the definition
of polygraphs:
\begin{lemma}\label{lemma:pushout}
The diagram 
\begin{displaymath}
\begin{xy}
\xymatrix{\sum_{S_n}\dglobst{n}\ar[r]^{\rho}\ar[d]_{\sum_{S_n}\cobnd{n}} & 
           \extend{n-1}{(\free S)}\ar[d]^{j_n}\\
          \sum_{S_n}\globst{n}\ar[r] & \extend{n}{(\free S)}}
\end{xy}
\end{displaymath}
is a pushout in $\compl$.
\end{lemma}

\section{Two classes of morphisms}\label{sec:twoclasses}

Let $\ctg{C}$ be a category, and $f:A\to B$, $g:C\to D$ morphisms. 
$f$ is {\em left-orthogonal} to $g$ (or, equivalently, 
$g$ is {\em right-orthogonal} to $f$) if, for each pair of morphisms
$u:A\to C$, $v:B\to D$ such that $g\circ u=v\circ f$, there exists
an $h:B\to C$ making the following diagram commutative:
\begin{displaymath}
\begin{xy}
\xymatrix{ A\ar[r]^u\ar[d]_f & C\ar[d]^g \\
           B\ar[r]_v\ar[ur]|h & D }
\end{xy}
\end{displaymath}

For any class $\mclass{M}$ of morphisms in $\ctg{C}$,
$\lortho{\mclass{M}}$ (resp.\ $\rortho{\mclass{M}}$) denotes
the class of morphisms in $\ctg{C}$ which are left-(resp.\ right-)
othogonal to all morphisms in $\mclass{M}$.

\subsection{Trivial fibrations}\label{subsec:fibrat}

Let $\cofib$ be the set $\set{\cobnd{n}| n\in\NN}$ as morphisms in $\compl$.
\begin{definition}
  A morphism of complexes is a {\em trivial fibration} if
  and only it belongs to $\rortho{\cofib}$.
\end{definition}
 In other words, $p:C\to D$ is a trivial fibration if for all $n$, $f:\dglobst{n}\to C$,
and $g:\globst{n}\to D$ such that $p\circ f=g\circ\cobnd{n}$, there is an $h:\globst{n}\to C$ making the following diagram commutative:
\begin{displaymath}
  \begin{xy}
    \xymatrix{\dglobst{n}\ar[r]^{f}\ar[d]_{\cobnd{n}} & C\ar[d]^{p}\\
              \globst{n}\ar[r]_{g}\ar[ur]|h & D}
  \end{xy}
\end{displaymath}

\begin{definition}\label{def:resolution}
Let $C$ be a complex. A {\em polygraphic resolution} of $C$ is
a pair $\pair{S}{p}$ where $S$ is a polygraph and
$p:\free S\to C$ is a trivial fibration.
\end{definition}
It was shown in~\cite{metayer:respol} that,
for each complex $C$, the counit of the
adjunction $\free{(.)} \dashv \abs{.}$,
\begin{displaymath}
\epsilon_C:\free{\abs{C}}\to C,
\end{displaymath} 
is a trivial fibration. Hence $\pair{\abs{C}}{\epsilon_C}$
is a polygraphic resolution of $C$, and we get the following result, which
will play an essential part in section~\ref{sec:main} below :

\begin{proposition}\label{prop:resolution}
Each complex $C$ has a polygraphic resolution.
\end{proposition}

\subsection{Cofibrations}\label{subsec:cofib}

\begin{definition}
A morphism of complexes is a {\em cofibration} if and
only if it is left-orthogonal to all trivial fibrations.
\end{definition}
Hence the class of cofibrations is exactly $\lortho{(\rortho{\cofib})}$.
Immediate examples of cofibrations are the $\cobnd{n}$'s themselves.
The following lemma summarizes standard properties of maps
defined by left-orthogonality conditions.
\begin{lemma}\label{lemma:leftlifting}
Let $\ctg{C}$ be a category, and $\mclass{M}$ an arbitrary class
of morphisms of $\ctg{C}$. Let $\mclass{L}=\lortho{\mclass{M}}$. Then
\begin{itemize}
\item $\mclass{L}$ is stable by direct sums: if $f_i:X_i\to Y_i$, $i\in I$
      is a family of maps of $\mclass{L}$ with direct sum
      $f=\sum_{i\in I} f_i:\sum_{i\in I}X_i\to\sum_{i\in I}Y_i$, then
      $f\in \mclass{L}$;
\item $\mclass{L}$ is stable by pushout: whenever $f\in\mclass{L}$
and
\begin{displaymath}
\begin{xy}
\xymatrix{X\ar[r]\ar[d]_f & Z\ar[d]^g\\
          Y\ar[r] & T}
\end{xy}
\end{displaymath}
is a pushout square in $\ctg{C}$, then $g\in\mclass{L}$;
\item suppose
\begin{displaymath}
\begin{xy}
\xymatrix{X_0\ar[r]^{l_0} & \cdots \ar[r]^{l_{n-1}}& X_n
\ar[r]^{l_n} & \cdots}
\end{xy}
\end{displaymath}
is a sequence of maps $l_n\in\mclass{L}$, with colimit $(X,m_n:X_n\to X)$.
Then  $m_0:X_0\to X$ belongs to $\mclass{L}$.
\end{itemize} 
\end{lemma}

\begin{proof}
We leave the first two claims as exercises. As for the third point,
let $f:Y\to Z$ be a morphism in $\mclass{M}$, and $u:X_0\to Y$,
$v:X\to Z$ such that the following diagram commutes:
\begin{displaymath}
\begin{xy}
\xymatrix{X_0\ar[d]_{m_0}\ar[r]^u & Y\ar[d]^f \\
          X\ar[r]_v   & Z}
\end{xy}
\end{displaymath}
Let us define $v_n=v\circ m_n$ for each $n\geq 0$. Thus,
for each $n\geq 0$,
$v_{n+1}\circ l_n=v\circ m_{n+1}\circ l_n=v\circ m_n=v_n$,
so that $v_n:X_n\to Z$
determines an inductive cone on the base $(X_n)$ to the vertex $Z$.
Let us define     
a family of maps $u_n:X_n\to Y$ satisfying
the following equations:
\begin{eqnarray}
f\circ u_n & = & v_n,\label{eq:fuv}\\
u_{n+1}\circ l_{n} & = & u_n\label{eq:ulu}.
\end{eqnarray}
Let $n=0$. Define $u_0=u$. We get
$f\circ u_0=f\circ u=v\circ m_0=v_0$, and~(\ref{eq:fuv})
holds. Thus $f\circ u_0=v_1\circ l_0$, and because $f\in\mclass{M}$ and
$l_0\in\mclass{L}$, there is an $u_1:X_1\to Y$ such that
$u_1\circ l_0=u_0$, so that~(\ref{eq:ulu}) holds. 
Suppose now that~(\ref{eq:fuv}) and~(\ref{eq:ulu}) hold for an $n\geq 0$.
By the induction hypothesis, the following diagram commutes:
\begin{displaymath}
\begin{xy}
\xymatrix{X_n\ar[r]^{u_n}\ar[d]_{l_n}& Y\ar[d]^f \\
          X_{n+1}\ar[r]_{v_{n+1}} & Z}
\end{xy}
\end{displaymath}          
with $f\in\mclass{M}$ and $l_n\in\mclass{L}$. Hence there is
a $u_{n+1}:X_{n+1}\to Y$ such that $f\circ u_{n+1}=v_{n+1}$
and $u_{n+1}\circ l_n=u_n$, and our equations hold for $n+1$.
In particular,~(\ref{eq:ulu}) means that $(u_n)$ determines
an inductive cone on the base $(X_n)$ to the vertex $Y$. As $X$
is the colimit of the $X_n$'s, there is a morphism $h:X\to Y$
such that, for each $n\geq 0$, $h\circ m_n=u_n$. In particular,
$h\circ m_0=u_0=u$.
Also, for each $n\geq 0$, $f\circ h\circ m_n=f\circ u_n=v_n=v\circ m_n$.
Uniqueness of connecting morphisms show that $f\circ h= v$.
Hence the following diagram is commutative
\begin{displaymath}
\begin{xy}
\xymatrix{X_0\ar[d]_{m_0}\ar[r]^u & Y\ar[d]^f \\
         X\ar[r]_v\ar[ur]|h  & Z}
\end{xy}
\end{displaymath}
and we have shown that $m_0\in\mclass{L}$, as required.
\end{proof}

\begin{definition}
A complex $C$ is {\em cofibrant} if $\initial\to C$
is a cofibration.
\end{definition}         

\begin{proposition}\label{prop:freecofib}
Free complexes are cofibrant.
\end{proposition}
\begin{proof}
Let $S$ be a polygraph and $C=\free S$. By lemma~\ref{lemma:pushout},
for each $n\geq -1$, the canonical inclusion $j_n:\extend{n}{C}\to
\extend{n+1}{C}$ is a pushout of $\sum_{S_n}\cobnd{n}$.
Now lemma~\ref{lemma:leftlifting} applies in the particular case
where $\mclass{L}$ is the class of cofibrations:
by the first point, 
$\sum_{S_n}\cobnd{n}$ is a cofibration, and by the second
point, so is $j_n$. By lemma~\ref{lemma:colimit}, $C$ is
a colimit of the sequence
\begin{displaymath}
\begin{xy}
\xymatrix{\extend{-1}{C}\ar[r]^{j_{-1}} & \extend{0}{C}
\ar[r]^{j_0} & \cdots\ar[r]^{j_{n-1}} & \extend{n}{C}\ar[r]^{j_n} & \cdots}
\end{xy}
\end{displaymath}
hence the third point of lemma~\ref{lemma:leftlifting} applies,
with $X_n=\extend{n-1}{C}$ and $l_n=j_{n-1}$, so that 
$\initial\to C$ is a cofibration. In other words, $C$ is cofibrant.
\end{proof} 

\section{Main result}\label{sec:main}

The main goal of this work is to establish the converse
of proposition~\ref{prop:freecofib}:

\begin{theorem}\label{thm:main}
Any cofibrant complex is isomorphic to a free one.
\end{theorem}

Let $\fcompl$ denote the full subcategory
of $\compl$ whose objects are the free complexes $\free S$
generated by polygraphs. Then, theorem~\ref{thm:main} reduces
to the following statement:

\begin{theorem}\label{thm:cauchy}
$\fcompl$ is Cauchy-complete.
\end{theorem}

Recall that a category $\ctg{C}$ is {\em Cauchy-complete} if 
all its idempotents split, that is, for each object $C$, and
each endomorphism $h:C\to C$ such that $h\circ h=h$, there
is an object $D$, together with morphisms $r:C\to D$, $u:D\to C$,
satisfying
\begin{eqnarray*}
u\circ r & = & h, \\
r\circ u & = & \id{}.
\end{eqnarray*}

Theorem~\ref{thm:cauchy} will be proved in annex~\ref{annex:cauchy}.
Let us assume the result for the moment, and let $C$ be a cofibrant
complex. By proposition~\ref{prop:resolution}, $C$ has a free resolution
$p:\free S\to C$, with $\free S$ an object of $\fcompl$. Because
$C$ is cofibrant, and $p$ is a trivial fibration, the identity
morphism $\id{C}:C\to C$ lifts through $p$, whence a morphism
$q:C\to \free{S}$ such that $p\circ q = \id{C}$. Let $h = q\circ p$,
$h\circ h= q\circ p\circ q\circ p=q\circ \id{C}\circ p=q\circ p=h$,
hence $h$ is an idempotent endomorphism of $\free S$. By using
Cauchy completeness, we get a polygraph $T$, and morphisms
$r:\free S\to\free T$, $u:\free T\to\free S$ such that $r\circ u=\id{\free{T}}$
and $u\circ r=h$. Now, let $f=p\circ u:\free T\to C$ and 
$g=r\circ q:C\to \free T$. We get
\begin{eqnarray*}
g\circ f & = & r\circ q\circ p\circ u\\
         & = & r\circ h \circ u\\
         & = & r\circ u\circ r\circ u\\
         & = & \id{\free{T}}\circ \id{\free{T}}\\
         & = & \id{\free{T}} .    
\end{eqnarray*}
Likewise 
\begin{eqnarray*}
f\circ g & = & p\circ u\circ r\circ q\\
         & = & p\circ h \circ q\\
         & = & p\circ q\circ p\circ q\\
         & = & \id{C}\circ \id{C}\\
         & = & \id{C}.    
\end{eqnarray*}
Hence $f:\free T\to C$ is an isomorphism with inverse $g=f^{-1}$
so that $C$ is isomorphic to a free object, as required.

\appendix

\section{Contexts}\label{annex:contexts}

\subsection{Indeterminates}\label{subsec:indeterminates}
Let $C$ be a complex, and $n\geq 1$. An {\em $n$-type} is an
ordered pair $\pair{x}{y}$ of parallel cells in $C_{n-1}$, that
is an element of $\paracell{n-1}{C}$. By~(\ref{eq:cpairs}),
$n$-types amount to morphisms $\theta:\dglobst{n}\to C$.
We shall use the same notations for
both sides of the natural equivalences~(\ref{eq:cyoneda})
and~(\ref{eq:cpairs}). 

\begin{definition}
The {\em type} of an $n$-cell $x\in C_n$ is
the pair $\pair{\sce{n-1}x}{\tge{n-1}x}$.
\end{definition}
Hence the type of an $n$-cell is a particular
$n$-type.

Given an $n$-type $\theta$, we may adjoin to $C$
an {\em indeterminate $n$-cell of type $\theta$} by
taking the following pushout in $\compl$:
\begin{displaymath}
\begin{xy}
\xymatrix{\dglobst{n}\ar[r]^{\theta}\ar[d]_{\cobnd{n}} & 
          C\ar[d]^{j_{\theta}} \\
         \globst{n}\ar[r]_{\vcx{x}} & C\ctxt{\vcx{x}}}
\end{xy}
\end{displaymath}
We let boldface variables $\vcx{x},\vcx{y},\ldots$ range
over indeterminates.

Let $\vcx{x}$ be an indeterminate $n$-cell of type $\theta$ and
$z:\globst{n}\to C$ an $n$-cell in $C$. To say that $z$ is
of type $\theta$ amounts to the commutativity
of the following diagram:
\begin{displaymath}
\begin{xy}
\xymatrix{\dglobst{n}\ar[r]^{\theta}\ar[d]_{\cobnd{n}} & C\ar[d]^{\id{}}\\
         \globst{n}\ar[r]_{z}& C}
\end{xy}
\end{displaymath}
The pushout property gives a unique morphism 
$\subst{z}:C\ctxt{\vcx{x}}\to C$ such
that $\subst{z}\circ\vcx{x}=z$ and $\subst{z}\circ j_{\theta}=\id{}$.
$\subst{z}$ is nothing but the operation of substituting
the cell $z$ for $\vcx{x}$ (see figure~\ref{fig:sub}).
\begin{figure}[ht]
\centering
\begin{displaymath}
\begin{xy}
\xymatrix{\dglobst{n}\ar[r]^{\theta}\ar[d]_{\cobnd{n}} & 
          C\ar[d]^{j_{\theta}}\ar[rdd]^{\id{}}&  \\
         \globst{n}\ar[r]^{\vcx{x}}\ar[rrd]_z&
          C\ctxt{\vcx{x}}\ar[rd]|{\subst{z}} &\\
         & & C}
\end{xy}
\end{displaymath}
\caption{substitution}
\label{fig:sub}
\end{figure}
Now $n$-cells of $C\ctxt{\vcx{x}}$ are formal composites of 
elements in $C_n\cup\set{\vcx{x}}$. Different expressions
may denote the same cell: however all those expressions 
contain the same number of occurrences of $\vcx{x}$. 
\begin{definition}
An {\em $n$-context over $\vcx{x}$} is an $n$-cell
of $C\ctxt{\vcx{x}}$ having exactly one occurrence
of $\vcx{x}$.
\end{definition}
We denote $n$-contexts over $\vcx{x}$ by
$c\ctxt{\vcx{x}},d\ctxt{\vcx{x}},\ldots$. 
An $n$-cell $z$ of $C$ is {\em adapted} to the context
$c\ctxt{\vcx{x}}$ if it has the same type as $\vcx{x}$.
Contexts are subject to the following operations:
\begin{itemize}
\item for each $n$-context 
      $c\ctxt{\vcx{x}}:\globst{n}\to C\ctxt{\vcx{x}}$ and each 
      adapted $n$-cell $z$, we denote by $c\ctxt{z}$
      the new $n$-cell of $C$
      obtained by substituting $z$ for $\vcx{x}$, in other terms
      $c\ctxt{z}=\subst{z}\circ c\ctxt{\vcx{x}}$;
\item let $u:C\to D$ be a morphism of complexes and 
      $c\ctxt{\vcx{x}}:\globst{n}\to C\ctxt{\vcx{x}}$ an $n$-context
      of $C$. Define a new indeterminate
      $\vcx{y}:\globst{n}\to D\ctxt{\vcx{y}}$ 
      by the following pushout square
       \begin{displaymath}
              \begin{xy}
                  \xymatrix{\dglobst{n}\ar[r]^{u\circ\theta}
                     \ar[d]_{\cobnd{n}} & 
                    D\ar[d]^{j_{u\circ\theta}} \\
                   \globst{n}\ar[r]_{\vcx{y}} & D\ctxt{\vcx{y}}}
              \end{xy}.
       \end{displaymath}
       This determines a unique morphism
       \begin{displaymath}
         \hat{u}:C\ctxt{\vcx{x}}\to D\ctxt{\vcx{y}}
       \end{displaymath}
       such that $j_{u\circ\theta}\circ u=\hat{u}\circ j_{\theta}$
       and $\hat{u}\circ\vcx{x}=\vcx{y}$ (see figure~\ref{fig:ctxtrans}).
      \begin{figure}[ht]
       \begin{displaymath}
        \centering
         \begin{xy}
          \xymatrix{\dglobst{n}\ar[r]^{\theta}\ar[d]_{\cobnd{n}} & 
          C\ar[d]^{j_{\theta}}\ar[rd]^{u} &  \\
         \globst{n}\ar[r]^{\vcx{x}}\ar[rrd]_{\vcx{y}}&
          C\ctxt{\vcx{x}}\ar[rd]|{\hat{u}} & D\ar[d]^{j_{\theta\circ u}}\\
         & & D\ctxt{\vcx{y}}}
        \end{xy}
        \end{displaymath}
        \caption{context transformation}
        \label{fig:ctxtrans}
        \end{figure}
       Thus $\hat{u}\circ c\ctxt{\vcx{x}}$ is an $n$-context
       over $\vcx{y}$, denoted by $c^u\ctxt{\vcx{y}}$.
\end{itemize}
Note that, if $z$ is an $n$-cell adapted to $c\ctxt{\vcx{x}}$, then
$u(z)$ is adapted to $c^u\ctxt{\vcx{y}}$ and
\begin{equation}
u(c\ctxt{z})=c^u\ctxt{u(z)}.\label{eq:context}
\end{equation}

\subsection{Thin contexts}\label{subsec:thin}

Let us introduce a few additional terminology about cells
and contexts. If $S$ is a polygraph, the elements of $\free S_n$
are the cells of {\em dimension $n$}, or $n$-cells. The 
{\em generators of dimension $n$}, or {\em $n$-generators} are
the elements of $S_n$. Each $n$-generator $\alpha$ determines
an $n$-cell $\free\alpha$. Such cells are called {\em atomic}.
All $0$-cells are atomic, and if $n>0$, each $n$-cell may
be expressed as a composition of atomic
cells and units on $n{-}1$-cells. For each $n$-cell $x$, and generator
$\alpha$, the number of occurrences of $\free\alpha$ in an expression
of $x$ only depends on $x$, not on the particular expression.
We call this number the {\em weight of 
$x$ at $\alpha$}, and denote it by $\wht{\alpha}{x}$.
The {\em total weight} of $x$ is
\begin{displaymath}
  \whtt{x} = \sum_{\alpha\in S_n}\wht{\alpha}{x}.
\end{displaymath}
The same definitions hold for contexts, where we take into account
all generators but the indeterminate. Thus, for instance.
$\whtt{\vcx{x}}=0$ for any indeterminate $\vcx{x}$. 

\begin{definition}\label{def:thin}
An $n$-context $c\ctxt{\vcx{x}}$ is {\em thin} if 
its total weight is zero.
\end{definition}
Now, if $x\in\free S_n$, either $\whtt{x}>0$ or there is a cell
$y\in\free S_{n{-}1}$ such that $x=\unit{n}(y)$. More generally,
if $\whtt{x}=0$, there is a unique integer $p<n$ with 
the following property:
\begin{itemize}
  \item there is a $p$-cell $z$ in $\free S$ such that
        $\whtt{z}>0$ and $x=\UNIT{n}{p}(z)$
\end{itemize}
Let us call $p$ the {\em thickness of $x$}, and denote it by
$p=\thk{x}$.
If $\whtt{x}\neq 0$, we define $\thk{x}=n$. The same
definitions immediately apply to contexts. In particular,
an $n$-context $c\ctxt{\vcx{x}}$ is thin if and only if 
$\thk{c\ctxt{\vcx{x}}}<n$. We finally associate to
each cell $x$ an integer $\size{x}$ by:
\begin{itemize}
  \item if $\whtt{x}\neq 0$, $\size{x}=\whtt{x}$;
  \item if $\whtt{x}=0$, $p=\thk{x}$, and $z$ is the unique
        cell in $\free S_p$ such that $x=\UNIT{n}{p}(z)$, then
        $\size{x}=\whtt{z}$.
\end{itemize}
In other words, the size of a cell $x$ is the number of generators of maximal
dimension needed to express $x$.
The size of contexts is defined accordingly.
Thus, the only contexts of size zero are just indeterminates.
We call those contexts {\em trivial}.

\begin{lemma}\label{lemma:scecontext}
  If $n>1$ and $c\ctxt{\vcx{x}}$ is a thin $n$-context,
  there is an $n{-}1$-context
  $d\ctxt{\vcx{y}}$ such that
  $\size{d\ctxt{\vcx{y}}}=\size{c\ctxt{\vcx{x}}}$ and
  for each adapted $n$-cell $z$,
\begin{displaymath}
          \sce{n-1}(c\ctxt{z}) = d\ctxt{\sce{n-1}(z)}.
\end{displaymath}
\end{lemma}
\begin{proof}
Let $\vcx{x}$ be an $n$-indeterminate of type
$\theta=\pair{x}{y}$.
We define a family $(\mclass{C}_i)_{0\leq i\leq n}$ 
of sets of $n$-contexts over $\vcx{x}$ by:
\begin{itemize}
  \item $\mclass{C}_0=\set{\vcx{x}}$;
  \item $\mclass{C}_{i+1}=\set{a\comp{i} c\ctxt{\vcx{x}}\comp{i}b}\cup
            \set{a\comp{i} c\ctxt{\vcx{x}}}\cup
            \set{c\ctxt{\vcx{x}}\comp{i}b}$, where
        $c\ctxt{\vcx{x}}\in\mclass{C}_i$, and $a$, $b$ are $n$-cells
        of $\free S$
        such that $a\seq i c\ctxt{\vcx{x}}\seq i b$, $\thk{a}>i$
        and $\thk{b}>i$.
\end{itemize}
Note that whenever $\thk{a}\leq i$ (resp.\ $\thk{b}\leq i$),
$a\comp i c\ctxt{\vcx{x}}=c\ctxt{\vcx{x}}$ (resp.\ $c\ctxt{\vcx{x}} \comp i b=c\ctxt{\vcx{x}}$).
Also the exchange rule allows to perform compositions
along higher dimensions outside those along lower dimensions. Hence
$\bigcup_{0\leq i\leq n}\mclass{C}_i$ contains all $n$-contexts
on $\vcx{x}$. As contexts in $\mclass{C}_n$ cannot be thin, all thin
contexts belong to $\bigcup_{0\leq i\leq n-1}\mclass{C}_i$.
Thus the lemma reduces to the following statement:
\begin{itemize}
  \item given $n>1$, $i\in\set{0,\ldots,n-1}$, and a thin
        $n$-context $c\ctxt{\vcx{x}}\in\mclass{C}_i$,
        there is an $n{-}1$-context $d\ctxt{\vcx{y}}$
        such that $\size{d\ctxt{\vcx{y}}}=\size{c\ctxt{\vcx{x}}}$ and,
        for each adapted $n$-cell $z$,
        $\sce{n-1}(c\ctxt{z}) = d\ctxt{\sce{n-1}(z)}$.
\end{itemize}
We prove this by induction on $i\in\set{0,\ldots,n-1}$.
\begin{itemize}
  \item If $i=0$, $c\ctxt{\vcx{x}}=\vcx{x}$ and we take
           $d\ctxt{\vcx{y}}=\vcx{y}$ where 
           $\vcx{y}$ is an $n{-}1$-indeterminate of type
           $\phi=\pair{\sce{n-2}(x)}{\tge{n-2}(x)}$.
           $\size{d\ctxt{\vcx{y}}}=\size{c\ctxt{\vcx{x}}}=0$ and 
           $\sce{n-1}(c\ctxt{z}) =\sce{n-1}(z)= d\ctxt{\sce{n-1}(z)}$.
  \item Suppose $0<i\leq n-1$ and the property holds for $i-1$.
        Let $c\ctxt{\vcx{x}}\in\mclass{C}_i$ be a thin context. 
        Then $c\ctxt{\vcx{x}}$ is of the form
        $a\comp{i-1} c'\ctxt{\vcx{x}} \comp{i-1} b$
        or  $a\comp{i-1} c'\ctxt{\vcx{x}}$
        or  $c'\ctxt{\vcx{x}} \comp{i-1} b$, 
        where $c'\ctxt{\vcx{x}}\in\mclass{C}_{i-1}$.
        We only treat the first case, the other two being very similar.
        Because $c\ctxt{\vcx{x}}$ is thin, so is $c'\ctxt{\vcx{x}}$,
        and $\whtt{a}=\whtt{b}=0$. Hence
        \begin{eqnarray*}
          \size{a} & = & \size{\sce{n-1}(a)}, \\
          \size{b} & = & \size{\sce{n-1}(b)}.
        \end{eqnarray*}
        By the induction hypothesis, there is an $n{-}1$-context
        $d'\ctxt{\vcx{y}}$ such that 
        $\size{d'\ctxt{\vcx{y}}}=\size{c'\ctxt{\vcx{x}}}$ and,
        for each adapted $n$-cell $z$,
        $\sce{n-1}(c'\tctxt{z}{\theta})=d'\tctxt{\sce{n-1}(z)}{\phi}$.
        As $i-1<n-1$, we may define
        \begin{displaymath}
          d\ctxt{\vcx{y}}= 
         \sce{n-1}(a)\comp{i}d'\ctxt{\vcx{y}}\comp{i}\sce{n-1}(b).
        \end{displaymath}
         We get
         \begin{eqnarray*}
           \size{d\ctxt{\vcx{y}}} & =
                   & \size{\sce{n-1}(a)}+\size{d'\ctxt{\vcx{y}}}+
                           \size{\sce{n-1}(b)}\\
                   & = & \size{a}+\size{d'\ctxt{\vcx{y}}}+\size{b}\\
                   & = & \size{a}+\size{c'\ctxt{\vcx{x}}}+\size{b}\\
                   & = & \size{c\ctxt{\vcx{x}}}
         \end{eqnarray*}
        and we get, for each adapted $n$-cell $z$,
        \begin{displaymath}
          \sce{n-1}(c\ctxt{z}) = d\ctxt{\sce{n-1}(z)}.
        \end{displaymath}
       \end{itemize}
\end{proof}

\begin{lemma}\label{lemma:trivial}
Let $c\ctxt{\vcx{x}}$ be an $n$-context and $z$ an adapted $n$-cell.
If $c\ctxt{z}=z$, then $c\ctxt{\vcx{x}}$ is trivial.
\end{lemma}

\begin{proof}
  By induction on the dimension $n$. If $n=1$,
  all contexts are trivial and we are done. Suppose now
  $n>1$ and the result holds in dimension $n-1$. Let
  $c\ctxt{\vcx{x}}$ be an $n$-context and $z$ an adapted $n$-cell
  such that
  \begin{equation}
    c\ctxt{z}=z \label{eq:czz}
  \end{equation}
  If $\whtt{c\ctxt{\vcx{x}}}>0$, then either $\whtt{z}>0$ and
  $\size{c\ctxt{z}}>\size{z}$, or $\whtt{z}=0$ and
  $\thk{c\ctxt{z}}=n>\thk{z}$. In both cases,
  $c\ctxt{z}\neq z$, a contradiction, and we are done.
  Otherwise, $c\ctxt{\vcx{x}}$ is thin, and lemma~\ref{lemma:scecontext} gives
  an $n{-}1$-context $d\ctxt{\vcx{y}}$ having the same size as $c\ctxt{\vcx{x}}$
  and satisfying $\sce{n-1}(c\ctxt{z})=d\ctxt{\sce{n-1}(z)}$.
  Hence, taking the source on both sides of~(\ref{eq:czz}), we get
  \begin{displaymath}
    d\ctxt{\sce{n-1}(z)}=\sce{n-1}(z).
  \end{displaymath}
   Thus, by the induction hypothesis, $d\ctxt{\vcx{y}}$ is a trivial
   $n{-}1$-context, hence $\size{d\ctxt{\vcx{y}}}=0$. Therefore
   $\size{c\ctxt{\vcx{x}}}=0$ as well, and $c\ctxt{\vcx{x}}$ is trivial.
\end{proof}

\bigskip

The following technical lemma will be crucial in the proof of
theorem~\ref{thm:cauchy}:

\begin{lemma}\label{lemma:thin}
Let $c\ctxt{\vcx{x}}$ be a thin $n$-context, and $z$ an adapted $n$-cell.
If $c\ctxt{z}$ is parallel to $z$, then $c\ctxt{z}=z$.
\end{lemma}

\begin{proof}
  Suppose $c\ctxt{\vcx{x}}$ is a thin $n$-context, and $z$ is an adapted
  $n$-cell such that $c\ctxt{z}\para z$. If $n=1$,
  thin contexts are trivial and the result is immediate.
  Otherwise, $n>1$ and by lemma~\ref{lemma:scecontext},
  there is an $n{-}1$-context $d\ctxt{\vcx{y}}$ such that 
  $\size{d\ctxt{\vcx{y}}}=\size{c\ctxt{\vcx{x}}}$ and
  $\sce{n-1}(c\ctxt{z})=d\ctxt{\sce{n-1}(z)}$. As $c\ctxt{z}$ is
  parallel to $z$, this implies $d\ctxt{\sce{n-1}(z)}=\sce{n-1}(z)$,
  and by lemma~\ref{lemma:trivial}, $d\ctxt{\vcx{y}}$ is a trivial context.
  Now $\size{d\ctxt{\vcx{y}}}=\size{c\ctxt{\vcx{x}}}=0$ so that
  $c\ctxt{\vcx{x}}$ is trivial,
  and $c\ctxt{z}=z$.
\end{proof}

\section{Cauchy completeness}\label{annex:cauchy}

This section is devoted to the proof of theorem~\ref{thm:cauchy}. 
Thus, let $S$ be a polygraph, and $h:\free S\to \free S$ an idempotent
morphism in $\compl$. We need to build a polygraph
$T$, together with morphisms $u:\free T\to \free S$ and
$r:\free S\to \free T$ such that
\begin{eqnarray}
  r\circ u & = & \id{},
   \label{eq:rui} \\
  u\circ r & = & h.
   \label{eq:urh}
\end{eqnarray}
We shall define $T$, $u$ and $r$ inductively on the dimension.
In dimension $0$, 
\begin{displaymath}
  T_0=\setbis{h(x)}{x\in \free S_0=S_0},
\end{displaymath}
$u$ is the inclusion $\free T_0=T_0\to \free S_0=S_0$, and for each
$x\in S_0$, $r(x)=h(x)$. The equations~(\ref{eq:rui}) and~(\ref{eq:urh})
are clearly satisfied. 

Suppose now that $n>0$ and $T$, $u$, $r$ have been defined up
to dimension $n{-}1$, and satisfy the required conditions. 
We shall extend the $n{-}1$ polygraph $T$ to an $n$-polygraph,
and the morphisms $u$, $r$ of $n{-}1$-complexes to morphisms
of $n$-complexes still satisfying the above equations.

\bigskip

\step{1} Let us split $S_n$ in three
subsets $S_n^0$, $S_n^1$ and $S_n^2$, according to the value of
$h(\free\alpha)$, for $\alpha\in S_n$:
\begin{itemize}
  \item $S_n^0=\setbis{\alpha\in S_n}{\whtt{h(\free\alpha)}=0}$, hence
        $S_n^0$ contains the generators $\alpha$ such that 
        $h(\free\alpha)$ is degenerate;
  \item $S_n^1$ contains the generators $\alpha\in S_n$ such that
        $\wht{\alpha}{h(\free\alpha)}=1$ and 
        $\wht{\beta}{h(\free\alpha)}=0$ if $\beta\notin S_n^0\cup\set{\alpha}$;
  \item $S_n^2=S_n\setminus{S_n^0\cup S_n^1}$. 
\end{itemize}
We may now define a set $T_n$ by:
\begin{displaymath}
  T_n=\setbis{h(\free\alpha)}{\alpha\in S_n^1}
\end{displaymath}
By definition, we get an inclusion map
\begin{displaymath}
  \upsilon: T_n\to \free S_n.
\end{displaymath}
such that
\begin{equation}
  h\circ \upsilon = \upsilon.
  \label{eq:hupsilon}
\end{equation}
Indeed, elements of $T_n$ belong to the image of the idempotent $h$, hence
are fixed by $h$.

We now define a graph $\sigma^T,\tau^T:\free T_{n{-}1}\OT T_n$ by
\begin{eqnarray}
  \sigma^T & = & r\circ \sce{n-1} \circ \upsilon 
   \label{eq:srsu}\\
  \tau^T & = & r\circ \tge{n-1} \circ \upsilon 
   \label{eq:trtu}  
\end{eqnarray}
where $\sce{n-1}$, $\tge{n-1}$ are the source and target maps in
$\free S$ and $r$ is given by the induction hypothesis:
\begin{displaymath}
  \begin{xy}
    \xymatrix{\free T_{n-1} & & T_n\ar[ll]_{\sigma^T,\tau^T}\ar[d]^{\upsilon}\\
       \free S_{n{-}1}\ar[u]^r & & \free S_n\ar[ll]^{\sce{n-1},\tge{n-1}}}
  \end{xy}.
\end{displaymath}
By using the fact that $r$ is a morphism up to dimension
$n{-}1$, we see that for each $\theta\in T_n$,
$\sigma^T(\theta)\para\tau^T(\theta)$ and the boundary conditions are
satisfied. Thus $T$ extends to an $n$-polygraph and the free
$n{-}1$-complex $\free T$ extends to a free $n$-complex.
We still denote these extensions by $T$, $\free T$, and
the source and target maps $\free{T}_{n-1}\OT\free{T}_n$ 
by $\sigma^T$ and $\tau^T$.

On the other hand, the following diagram commutes
\begin{displaymath}
  \begin{xy}
    \xymatrix{\free T_{n-1}\ar[d]_u & T_n\ar[l]_{\sigma^T}\ar[d]^{\upsilon}\\
       \free S_{n{-}1} & \free S_n\ar[l]^{\sce{n-1}}}
  \end{xy}
\end{displaymath}
because  
\begin{eqnarray*}
  u\circ \sigma^T & = & u\circ r\circ \sce{n-1}\circ\upsilon,\\
                & = & h\circ \sce{n-1}\circ\upsilon,\\
                & = & \sce{n-1}\circ h\circ \upsilon,\\
                & = & \sce{n-1}\circ\upsilon. 
\end{eqnarray*}
Likewise
\begin{eqnarray*}
  u\circ\tau^T & = & u\circ r\circ \tge{n-1}\circ\upsilon.
\end{eqnarray*}
Hence $\upsilon:T_n\to \free S_n$ gives rise to
$u_n:\free T_n\to \free S_n$, extending $u$ to a morphism
of $n$-complexes $\free T\to \free S$.

To sum up, we have extended $T$ and $u$ up to dimension $n$.
Remark that the only property of $T_n$ we needed so far is
that its elements are fixed by $h$.

\bigskip

\step{2} We introduce an auxiliary $n$-polygraph $U$ by
\begin{itemize}
  \item $U$ is identical to $S$ up to dimension $n{-}1$;
  \item $U_n=S_n^0+S_n^1$ and the  source and target maps
        $\free U_{n-1}\OT U_n$ simply restrict those on $S_n$.
\end{itemize}
Thus we get an inclusion monomorphism  of
$n$-polygraphs $\iota:U\to S$, generating a monomorphism of
$n$-complexes $\free\iota:\free U\to\free S$. The restrictions
of $\sce{n-1}$ and $\tge{n-1}$ to $\free U_{n}$ will be denoted
by $\sigma^U$ and $\tau^U$, as well as the correponding maps
on generators: $\free U_{n-1}\OT U_n$.
\begin{lemma}\label{lemma:U}
  There are morphisms of $n$-complexes
  \begin{displaymath}
    h':\free U\to \free U,\quad
    k : \free S\to\free U,
  \end{displaymath}
 such that the following diagram commutes:
  \begin{displaymath}
    \begin{xy}
      \xymatrix{\free U\ar[r]^{\free{\iota}}\ar[d]_{h'} &
         \free S\ar[dl]|k\ar[d]^h \\
                \free U\ar[r]_{\free{\iota}} & \free S}
    \end{xy}
  \end{displaymath}
\end{lemma}
\begin{proof}
  The existence of $h'$ making the outer square commutative
  follows from the remark that $\free U$ is stable by $h$, 
  so that $h'$ is simply the restriction of $h$ to $\free U$.

  As for $k$, the statement reduces to the fact that all $n$-cells
  of the form $y=h(x)$ in $\free S$ can be expressed by generators
  taken from $U_n$.  Thus, let $y=h(x)$ an $n$-cell of $\free S$.
  The (occurrences) of generators $\gamma\in S_n$
  such that $\wht{\gamma}{y}>0$ may
  be arranged in a list
  \begin{displaymath}
    \alpha_1,\ldots,\alpha_p,\beta_1,\ldots,\beta_q
  \end{displaymath}
  where $\alpha_i\notin S_n^0$ and $\beta_j\in S_n^0$. Notice
  that repetitions are possible.
  Let $y_i=h(\free\alpha_i)$ for each $i\in\set{1,\ldots,p}$.
  As $h(y)=y$ and $\whtt{h(\free\beta_j)}=0$ for each
  $j\in\set{1,\ldots,q}$, we get
    \begin{equation}
      \sum_{\gamma\notin S_n^0}\wht{\gamma}{y}
      = \sum_{i=1}^{i=p}\sum_{\gamma\notin S_n^0}\wht{\gamma}{y_i}
     \label{eq:sumw}
    \end{equation}
 But $h(y_i)=y_i=h(\free\alpha_i)$ and the generators of $y_i$
 cannot be all in $S_n^0$, otherwise $\whtt{h(y_i)}=0$, in contradiction
 with $\alpha_i\notin S_n^0$. Hence, for each 
 $i\in\set{1,\ldots,p}$, there is at leat one $\gamma\notin S_n^0$,
 such that $\wht{\gamma}{y_i}>0$. Therefore, the left hand side
 of~(\ref{eq:sumw})
 is equal to $p$, whereas the right hand side has $p$ terms, all
 of which are $\geq 1$. This is possible only if
 \begin{displaymath}
   \sum_{\gamma\notin S_n^0}\wht{\gamma}{y_i} = 1
 \end{displaymath}
 for each $i\in\set{1,\ldots,p}$. Let $\delta_i$ be the only
 generator in $S_n\setminus S_n^0$ such that $\wht{\delta_i}{y_i}=1$.
 The occurrences of $n$-generators in $y_i$ are exactly those
 in $h(y_i)$, hence in $h(\free\delta_i)$. It follows that
 $\wht{\delta_i}{h(\free\delta_i)}=1$ and $\wht{\gamma}{h(\free\delta_i)}=0$
 for each $\gamma\notin S_n^0\cup\set{\delta_i}$. This means  exactly
 that $\delta_i\in S_n^1$. Therefore $y$ can be expressed by using 
 as $n$-generators
  \begin{displaymath}
    \delta_1,\ldots,\delta_p,\beta_1,\ldots,\beta_q,
  \end{displaymath} 
 all in $S_n^0\cup S_n^1=U_n$. Thus for each $x\in\free S_n$,
 there is a unique $y\in\free U_n$ such that $\free\iota(y)=h(x)$.
 Hence a morphism $k:\free S\to \free U$
 such that $\free\iota\circ k=h$. Finally 
 $\free\iota\circ k\circ \free\iota= h\circ\free\iota = \free\iota\circ h'$,
 and because $\free\iota$ is a monomorphism, $k\circ\free\iota=h'$. 
\end{proof}

\bigskip

If $x\in \free T_n$, $u(x)\in \free S_n$ can be expressed by
generators from $U_n$, hence there is a morphism
$u':\free T\to \free U$ such that $u=\free\iota\circ u'$. Of course
$u'$ coincides with $u$ in all dimensions $i<n$.

\bigskip

\step{3} We now define a morphism $r':\free U\to\free T$ which coincides
with $r$ in dimensions $i<n$. All we need is a map
\begin{displaymath}
  \rho:U_n\to \free T_n
\end{displaymath}
satisfying the boudary conditions. Thus, let $\alpha\in U_n$, we distinguish
two cases, according as $\alpha\in S_n^0$ or $\alpha\in S_n^1$.

\case{1} Let $\alpha\in S_n^0$. There is a unique $y\in\free S_{n-1}$
such that $h(\free\alpha)=\unit{n}(y)$. Now $r(y)\in\free T_{n-1}$, so that 
we may define $\rho(\alpha)=\unit{n}(r(y))$. The boundary conditions
are straightforward in this case.

\case{2} Let $\alpha\in S_n^1$. There is a unique generator $\theta\in T_n$
such that $h(\free\alpha)=\upsilon(\theta)$. We define
$\rho(\alpha)=\free\theta$. By using the induction hypothesis on $r$
and $u$, we get
\begin{eqnarray*}
  \sigma^T(\rho(\alpha)) & = & \sigma^T(\free\theta) \\
                       & = & r(\sce{n-1}(\upsilon(\theta))) \\
                       & = & r(\sce{n-1}(h(\free\alpha))) \\
                       & = & r(h(\sce{n-1}(\free\alpha))) \\
                       & = & r(u(r(\sce{n-1}(\free\alpha)))) \\
                       & = & r(\sce{n-1}(\free\alpha)) \\
                       & = & r'(\sigma^U(\alpha))
\end{eqnarray*}
Hence $\sigma^T(\rho(\alpha))=r'(\sigma^U(\alpha))$ and likewise
 $\tau^T(\rho(\alpha))=r'(\tau^U(\alpha))$, and the boundary
conditions are satisfied.

Thus $\rho$ gives rise to a morphism of complexes
$r':\free U\to\free t$ extending $r$ up to dimension $n$.

\bigskip

\step{4} Having defined $u':\free T\to \free U$ and $r':\free U\to \free T$,
we first note that $u'\circ r'=h'$, which directly follows from
our definition of $r'$. We now prove the following lemma:
\begin{lemma}\label{lemma:r'u'i}
$r'\circ u'=\id{}$.
\end{lemma}
\begin{proof}
 $r'\circ u'$ is an endomorphism of the complex $\free T$. We know
 by the induction hypothesis that $r'\circ u'=r\circ u=\id{}$ 
 in all dimensions $i<n$. Thus, it suffices to show that, for
 each generator $\theta\in T_n$,
 \begin{equation}
   r'(u'(\free\theta)) = \free\theta
    \label{eq:r'u'i}
 \end{equation}
This follows from two facts:
\begin{itemize}
  \item the two members of~(\ref{eq:r'u'i}) are parallel cells:
    \begin{displaymath}
      \sigma^T(r'(u'(\free\theta))) = r'(u'(\sigma^T(\free\theta)))
    \end{displaymath}
   because $r'$, $u'$ are morphisms. But
   $\sigma^T(\free\theta)$ has dimension $n-1$, where, by
   the induction hypothesis, $r'\circ u'=\id{}$,
   so that the above equation becomes
 \begin{displaymath}
      \sigma^T(r'(u'(\free\theta))) = \sigma^T(\free\theta)
    \end{displaymath}
  and likewise 
  \begin{displaymath}
      \tau^T(r'(u'(\free\theta))) = \tau^T(\free\theta).
    \end{displaymath}       
\item there is a {\em thin} $n$-context $c\ctxt{\vcx{x}}$ in $\free T$
      such that
      \begin{displaymath}
        r'(u'(\free\theta))=c\ctxt{\free\theta}.
      \end{displaymath}
  In fact, by definition of $T_n$, there is a generator $\alpha\in S_n^1$
  such that $u'(\free\theta)=h(\free\alpha)$. As $h(\free\alpha)$ contains
  a single occurrence of $\free\alpha$, there is an $n$-context 
  $d\ctxt{\vcx{y}}$ in $\free U$ such that
  $u'(\free\theta)=d\ctxt{\free\alpha}$. Now by applying~(\ref{eq:context})
  of section~\ref{subsec:indeterminates},
  \begin{eqnarray*}
    r'(d\ctxt{\free\alpha}) & = & d^{r'}\ctxt{r'(\free\alpha)}\\
                            & = & d^{r'}\ctxt{\rho(\alpha)} \\
                            & = & d^{r'}\ctxt{\free\theta}
  \end{eqnarray*}
  Define $c\ctxt{\vcx{x}}=d^{r'}\ctxt{\vcx{x}}$. All generators
  of $d\ctxt{\free\alpha}$ but $\alpha$ itself belong to $S_n^0$,
  hence are sent to identities by $r'$. Therefore $c\ctxt{\vcx{x}}$ is
  thin, and we are done.
  \end{itemize}
 $c\ctxt{\vcx{x}}$ is a thin context such that
 $c\ctxt{\free\theta}\para\free\theta$. By lemma~\ref{lemma:thin},
 $c\ctxt{\free\theta}=\free\theta$ and~(\ref{eq:r'u'i}) is proved.
 \end{proof}

\bigskip

\step{5} We complete the argument by defining $r=r'\circ k$. Hence
$r$ is a morphism $\free S\to \free T$ and
\begin{eqnarray*}
  u\circ r & = & \free\iota\circ u'\circ r'\circ k, \\
           & = & \free\iota\circ h'\circ k,\\
           & = & \free\iota\circ k\circ\free\iota\circ k,\\
           & = & h\circ h, \\
           & = & h.
\end{eqnarray*}
Also
\begin{eqnarray*}
  r\circ u & = & r'\circ k\circ \free\iota\circ u',\\
           & = & r'\circ h'\circ u',\\
           & = & r'\circ u'\circ r'\circ u',\\
           & = & \id{}\circ\id{},\\
           & = & \id.
\end{eqnarray*}
Thus~(\ref{eq:rui}) and~(\ref{eq:urh}) hold in dimension $n$ and we are done.
 \end{document}